\documentclass[12pt]{article}

\UseRawInputEncoding

\usepackage{mathtools,amsthm,amssymb,color,cases}
\usepackage{mathrsfs}
\usepackage[abbrev]{amsrefs}
\usepackage{tikz}
\usepackage{enumitem}
\usepackage{setspace}

\usepackage{latexsym}
\usepackage{mathrsfs}
\usepackage{comment}
\usepackage{microtype}

\usepackage{indentfirst}

 \newtheoremstyle{mystyle}% % Name
    {}%                      % Space above
    {}%                      % Space below
    {\normalfont}%           % Body font
    {}%                      % Indent amount
    {\bfseries}%             % Theorem head font
    {}%                      % Punctuation after theorem head
    { }%                     % Space after theorem head, ' ', or \newline
    {}%                      % Theorem head spec (can be left empty, meaning `normal')
  \theoremstyle{mystyle}     

\usepackage{geometry}
\allowdisplaybreaks[1]
\theoremstyle{definition}
\newtheorem{thm}{Theorem}[section]

\newtheorem{lem}[thm]{Lemma}
\newtheorem{rem}[thm]{Remark}
\newtheorem*{rem*}{Remark}

\newtheorem*{definition}{Definition}

\newtheorem{prop}[thm]{Proposition}

\newtheoremstyle{part}{}{}{\normalfont}{}{\itshape}{.}{.5em}{}
\theoremstyle{part}

\newtheorem{step}{Step}[section]

\newcommand\blfootnote[1]{%
  \begingroup
  \renewcommand\thefootnote{}\footnote{#1}%
  \addtocounter{footnote}{-1}%
  \endgroup
}

\numberwithin{equation}{section}
\numberwithin{thm}{section}

\newenvironment{equ}{
    \begin{equation}
}{
    \end{equation}
}
\newenvironment{equ*}{
    \begin{equation*}
}{
    \end{equation*}
}

\newtheoremstyle{part}{}{}{\normalfont}{}{\itshape}{.}{.5em}{}
\theoremstyle{part}

\usepackage{etoolbox}
\AtEndEnvironment{proof}{\setcounter{proofpart}{0}}

\newcommand{\pt}{\partial}

\newcommand{\F}{\mathcal {F}}

\newcommand{\Z}{\mathbb{Z}}

\newcommand{\R}{\mathbb R}

\newcommand{\intr}{\int_{\R^2}}      
\newcommand{\intt}{\int^t_0}
\newcommand{\nonlin}{\intt \P(t-s) \ast (u(s) \cdot \N) \t(s) \dd s}
\newcommand{\lin}{\P(t) \ast \t_0}

\newcommand{\al}{\alpha}

\newcommand{\x}{\xi}

\renewcommand{\t}{\theta}
\renewcommand{\th}{\hat\theta}

\newcommand{\N}{\nabla }

\def\<{\langle }

\newcommand{\supp}{\text{ supp }}

\renewcommand{\P}{P_{\al/2}}

\newcommand*{\dd}{\mathop{}\!\mathrm{d}}

\DeclarePairedDelimiter{\norm}{\lVert}{\rVert}
\DeclarePairedDelimiter{\bignorm}{\bigg\|}{\bigg\|}
\DeclarePairedDelimiter{\quot}{``}{"}

\DeclarePairedDelimiter{\Bigc}{\Big{(}}{\Big{)}}
\DeclarePairedDelimiter{\Biggc}{\Bigg{(}}{\Bigg{)}}

\DeclarePairedDelimiter{\Bigf}{\Big{[}}{\Big{]}}



\begin{document}

\begin{center}
    {\bf \large Large-Time Behaviour of Solutions to the Surface Quasi-Geostrophic Equations} \\
    
    \phantom{}
    
    D\'aith\'i \'O hAodha* \quad Tsukasa Iwabuchi** \\
    
     \phantom{}
    
    Mathematical Institue, Tohoku University,
    
    980-0845
\end{center}

\blfootnote{Email: *david.declan.hughes.p6@dc.tohoku.ac.jp, **t-iwabuchi@tohoku.ac.jp}

\begin{center}
\begin{minipage}{135mm}
\footnotesize
{\sc Abstract. }
We construct a linear approximation of the solution to the Surface Quasi-Geostrophic Equations in $\R^2$, and obtain a convergence rate in $L^p$ between the solution and this approximation with respect to time. We also demonstrate that the nonlinear term of the solution is bounded sharply in $L^p$ by the same function of time.
\end{minipage}
\end{center}

\section{Introduction}

This paper is concerned with the Surface Quasi-Geostrophic Equation.
\begin{equ} \label{(Q)} \begin{cases}
\partial_t \theta + (- \Delta)^{\alpha/2} \theta + (u\cdot \nabla) \theta = 0, 
& \ \text{in} \ (0,\infty) \times \mathbb{R}^2, \\
u = (-R_2\theta, R_1\theta),
& \ \text{in} \ (0,\infty) \times \mathbb{R}^2, \\
\theta|_{t=0} = \theta_0, 
& \ \text{in} \ \mathbb{R}^2. 
\end{cases} \vspace{-2pt}
\end{equ}
\\
Here, $\theta:  (0,\infty) \times \mathbb{R}^2 \rightarrow \mathbb{R}$ is an unknown function, representing the potential temperature of a fluid parcel at a point $(t,x)$ in spacetime; and $u$ represents the velocity of a fluid parcel. $R_j = \partial_j (-\Delta)^{-1/2} = \F^{-1}[ \, i\x_j / |\x| \, ]\ast $ is the $j$-th Riesz transform; and $\alpha \in [1,2]$. 
We refer to references~\cites{surface1, surface2, Pedlosky, landau, Const Crit, Blumenthal} for the physical meaning and derivation of the equations.
We will prove existence and uniqueness results for given initial data; and, under some slightly stronger restrictions, we will consider the large-time behaviour of solutions.

Let us recall several existing results related to the regularity of solutions. In the subcritical case, $\al \in (1,2],$ unique global existence and regularity can be shown, for initial data 
\[ \t_0 \in L^1 \cap L^p, \ p\in \Big{(}\frac{2}{\al-1}, \infty \Big{]},\]
by classical methods using the Banach fixed point theorem, similarly to \cites{Escobedo-Zuazua, Iwabuchi-2015}. Global unique existence and regularity for the subcritical case are also proven on the torus in \cite{Constantin}.

In the critical case, $\alpha = 1$, we can guarantee local existence of a solution for uniformly continuous 
initial data (see~\cite{WZ-2011}, which studies the problem
in the framework of Besov spaces larger than $L^\infty$). 
The local smooth solutions are then extended to global smooth solutions, 
as in the papers~\cites{Kiselev-Nazarov-Volverg,CaVa-2010,Vicol}.

Finally, for the supercritical case, $\al \in [0,1),$ the paper \cite{CotVic-2016} proves global regularity for all $\al \in [\al_0,1)$, where $\al_0$ grows with respect to the size of the initial data. Global regularity for large data in the supercritical case is an open problem. In this paper, we discuss a sharp decay estimate of
the nonlinear part of the solution in the subcritical and critical cases. We also construct a linear approximation of the solution in $L^p$.

We begin our study of \eqref{(Q)} by defining the following function
\begin{align} \label{heat}
P_{\al/2}(t,x) = \F^{-1}\left[e^{-t|\xi|^\alpha} \right](x), \text{ for } t>0, \ x\in\R^2,
\end{align}
which is the fundamental solution to the fractional heat equation, the linear part of \eqref{(Q)}. We will use this new function to introduce the idea of mild solutions.

For this paper, we set the initial data as follows
\begin{equation}\t_0 \in W^{1,1} \cap W^{1,\infty}, \label{data} \end{equation}
which provides us with sufficient regularity for global existence and smoothness of solutions in the subcritical and critical cases.

\begin{definition} (Mild Solution)
A function, $\t$, is a mild solution of \eqref{(Q)} if
\begin{align} 
\label{integral eqn}
& \t(t) = \lin - \nonlin, \text{ for all } t>0, \\
& \lim_{t\to 0^{+}} \t(t) = \t_0, \text{ in } L^p, \text{ for all } p\in[1,\infty), \\
& \theta \in C([0,\infty) ; L^p(\R^2)), \text{ for all } p \in [1,\infty), \\
& \theta \in C((0,\infty) ; W^{2,p}(\R^2) ) \cap C^1((0,\infty) ; L^p(\R^2)), \text{ for all } p \in [1,\infty].
\end{align}
\end{definition}

We also define $M$, the $\quot{\text{mass}}$ of the solution, $\t$, as 
\[
M := \intr \t_0(x) \dd x,
\]
and denote the linear part of $\t$ by
\[
U(t) := \P(t) \ast \t_0.
\]
We state an existence result for global solutions.
\begin{prop}
Let $\alpha \in [1,2]$ and $\theta _ 0 \in W^{1,1} \cap W^{1,\infty}$. 
Then there exists a unique global mild solution 
$\theta \in C( [0,\infty) ; W^{1,p}  ) 
\cap C((0,\infty) ; W^{1,1} \cap W^{1,\infty})$ of \eqref{(Q)}, 
for all $1 \leq p < \infty$. 
\end{prop}

As a first step, it is possible to prove that, for 
$1 \leq p \leq \infty$,
\[
\lim_{t\to\infty}t^{\frac{2}{\alpha}(1-\frac{1}{p})}
\| \theta (t) - M P_{\alpha /2} (t) \|_{L^p} = 0 . 
\]
This explains that the $L^p$-decay of the solution is 
essentially equivalent to that of the linear solution.
In this paper, we will extract the optimal decay of the nonlinear part of $\t$. 
Our result reads as follows:

\begin{thm} \label{main theorem}
Let $\al \in [1,2]$, $p \in [1,\infty]$ and 
\begin{equation}\label{b}
b_{\alpha ,p}(t) = 
\begin{cases}
t^{\frac{2}{\alpha}(1-\frac{1}{p})+\frac{3}{\alpha}-1} & \text{ if } \alpha \in (1,2], 
\\
t^{2(1-\frac{1}{p})+2} \ln t & \text{ if } \alpha = 1.  
\end{cases}
\end{equation}
Let $\t_0 \in W^{1,1} \cap W^{1,\infty}$, 
and let $\theta$ be a mild solution. Also assume $|x|^2 \t_0 \in L^1(\R^2) $. 
Then we have the following convergence:
\begin{align} \label{thm3 decay}
b_{\alpha ,p}(t) % t^{\frac{2}{\al} (1-\frac{1}{p}) + \frac{3}{\al} - 1} 
    & \bignorm{ \theta(t) -M\P(t) + \N\P(t) \cdot \intr y \t_0(y) \dd y - 
       \sum_{i,j = 1}^2 \partial_{x_i} \partial_{x_j} \P(t) \intr y_i y_j \t_0 \dd y  \\
    & + \intt \P(t-s) \ast \N \cdot \Big{(} (RU(s))U(s) \Big{)} \dd s}_p
    \to 0, \text{ as } t \to \infty. 
     \nonumber 
\end{align}
Furthermore, there exists $\theta_0$ such that the nonlinear component is optimally bounded by $b_{\alpha ,p}(t)$ in the $p=2$ case. That is 
\begin{align}  \label{bound from below}
    & \bignorm{\intt \P(t-s) \ast \N \cdot \Big{(} (RU(s))U(s) \Big{)} \dd s}_2 \simeq 
    \frac{1}{b_{\alpha ,2}(t)} 
\end{align}
for all sufficiently large $t$. 
\end{thm}

Our proof of \eqref{bound from below} is fully self-contained, with all necessary work shown explicitly in this paper.
For the bound from above in \eqref{bound from below}, the proof mainly consists of carefully taking Besov norms of the solution via Littlewood-Paley decomposition; and, in the case of the $L^p$ norm for $1\leq p <2$, applying the Hardy-Littlewood-Sobolev and Gr\"onwall inequalities.
For the bound from below, we force the initial data to take a shape resemblant of a Gaussian function with a narrow support, and exploit this assumption to achieve the necessary inequalities.

The first term inside the norm of \eqref{thm3 decay} is the solution, $\t$. 
The second, third, and fourth terms are a (Taylor-expanded) approximation of the linear term of $\t$.
The final term is a linear approximation of the nonlinear term in $\t$.

\begin{rem}
If the initial data, $\t_0$, is radially symmetric, the approximation of the nonlinear term becomes $0$ (see e.g. \cite{Giga-Giga-Saal}, page 46). That is,
\[
\intr \P(t-s) \ast (RU(s)\cdot\nabla)U(s) \dd s = 0.
\]
\end{rem}

\noindent 
{\bf Notation}. 
For a function, $f$, we denote the Fourier transform of $f$ as follows:
\[
\F[f](\x) := \hat{f}(x) := \frac{1}{2\pi}\intr e^{-i x\cdot\x} f(x) \dd x.
\]
The inverse Fourier transform is then written as 
\[
\F^{-1}[\hat{f}](x) := \frac{1}{2\pi} \intr e^{i x \cdot \xi} \hat{f}(\x) \dd \x.
\]
For the purpose of calculating inequalities, we will frequently omit the factor of $1/2\pi$, as it will have no influence on the proofs.
Let $\mathcal S' = \mathcal S'(\mathbb R^2)$ be the space of tempered distributions. 
Let $\mathcal P = \mathcal P (\mathbb R^2)$ be the set of all polynomials.

\section{Preliminaries}

We recall the definition and some basic properties of Besov spaces, 
and write the $L^p$-norm decay of the solution, $\t$.

\subsection{Besov Spaces}

We use the Littlewood-Paley decomposition of unity to define homogeneous Besov spaces. 

\begin{definition}

Let $\{ \phi_k \}_{k\in\Z}$ be a set of non-negative measurable functions such that 
\begin{enumerate}
    \item $\displaystyle \sum_{k\in\Z} \hat{\phi}_k (\x) = 1, \text{ for all } \x \in \R^2 \backslash \{0\}$,
    \item $\hat{\phi}_k(\x) = \hat{\phi}_0(2^{-k}\x)$,
    \item $\supp \hat{\phi}_k (\x) \subseteq \{ \x \in \R^2 \ | \ 2^{k-1} \leq |\x| \leq 2^{k+1} \}$.
\end{enumerate}
The Besov norm is then defined as follows. 
For $f \in \mathcal{S'}/\mathcal{P}$, $1 \leq p,q \leq \infty$, and $s \in \R$,
\[
\norm{f}_{\dot{B}^{s}_{p,q}} := \bignorm{\{ 2^{sk} \norm{\phi \ast f}_{p} \}_{k\in\Z}}_{l^q}.
\]
Finally, the set $\dot{B}^{s}_{p,q}$ is defined as the set of distributions, $f \in \mathcal{S'}/\mathcal{P}$, whose Besov norm is finite.

\end{definition}

We introduce the following propositions, and refer to \cite{Triebel} for their proofs.

\begin{prop} 

Let $1 \leq p,q \leq \infty$, and $s \in \R$. Then for $f \in \dot{B}^{s+1}_{p,q}$,
\[
\norm{\nabla f}_{\dot{B}^{s}_{p,q}} \leq C \norm{f}_{\dot{B}^{s+1}_{p,q}}.
\]

\end{prop}

\begin{prop} 

Let $1 \leq p \leq \infty$. Then for $f \in \dot{B}^{2(1-\frac{1}{p})}_{1,1}$,
\[
\norm{f}_{\dot{B}^{0}_{p,1}} \leq C \norm{f}_{\dot{B}^{2(1-\frac{1}{p})}_{1,1}}.
\]

\end{prop}

\subsection{$L^p$-Norm Decay}

\begin{prop}(Hardy-Littlewood-Sobolev~\cite{Hedberg}) \label{HLS}
Let $0 < \al < n$ and $1<p<r<\infty$, such that $\frac{1}{r} = \frac{1}{p} - \frac{\al}{n}$. Then there exists a constant $ C > 0$ such that
\[
\norm{ (-\Delta)^{-\al/2} f}_{L^r(\R^n)} \leq C \norm{f}_{L^p(\R^n)}.
\]
\end{prop}

\begin{prop}
Let $n \in \mathbb N$, $p \in [1,\infty]$, $k \in \mathbb N$, $\alpha \in [1,2]$. 
Then there exists $C>0$ such that, for all $1 \leq j \leq n$, and $t > 0  $, 
\begin{align} \label{P-estimate}
\norm{\pt_j^k \P(t)}_{L^p(\R^n)} \leq C t^{ \frac{n}{\al}(1-\frac{1}{p}) - \frac{k}{\al}}.
\end{align}
The decay rate is easily obtained by a change of variables, recalling \eqref{heat}. The overall boundedness  bis proven using the Hausdorff-Young inequality. 
\end{prop}

\begin{prop}\label{prop:decay} Let $\theta_0 \in W^{1,1} \cap W^{1,\infty}$, and $1 \leq p \leq \infty$. 

\begin{enumerate}[label=(\roman*)]

\item 
For all $\alpha \in [1,2]$, there exists $C >0$ such that, for all $t>0$,
\begin{align} \label{L^p decay}
    \norm{\t(t)}_{p} \leq C (  t + 1 )^{-\frac{2}{\al}(1-\frac{1}{p}) }.
\end{align}

\item
For $\al = 1$, $\beta > 0$, there exists $C_\beta > 0$ such that, for all $t \geq 1$,
\begin{align} \label{Iwab Lp deriv decay}
    \norm{ |\N|^{\beta} \t(t)}_p \leq C_{\beta} t^{-2(1-\frac{1}{p}) - \beta }.
\end{align}

\end{enumerate}
\end{prop}

For the proof of (i), see \cite{Cordoba}. For the proof of (ii), see Proposition 4.3 in \cite{Iwabuchi-2020}.

\section{Large-Time Behaviour}

We will now begin to discuss the large-time behaviour of the solution, $\t$. Before beginning our proof of Theorem \ref{main theorem}, we begin with a less strong approximation of the solution. Results similar to the below have been proven with respect to the $N$-dimensional convection-diffusion equations~\cite{Escobedo-Zuazua} and the critical Burger's equations~\cite{Iwabuchi-2015}.

\subsection{Approximation by the Fractional Heat Kernel}

\begin{prop}
Let $\al \in [1,2]$ and $p \in [1,\infty]$. 
Let $\t_0 \in W^{1,1} \cap W^{1,\infty}$,
and also assume $|x| \t_0 \in L^1$. 
Then the solution, $\t$, to (\ref{integral eqn}) satisfies
\begin{align} \label{thm2fast}
    \norm{ \t(t) - M\P(t) }_p \leq C t^{-\frac{2}{\al} (1-\frac{1}{p}) - \frac{1}{\al}} , \text{ for all } t \geq 1. 
\end{align}
\end{prop}

In order to prove the above proposition, it is useful to split the norm into linear and nonlinear parts as follows: 
\begin{align*}
    \norm{\t(t) - M\P(t)}_p \leq \norm{\P(t) \ast \t_0 - M\P(t)}_p + \bignorm{ \nonlin }_p . 
\end{align*}
We then prove the bound \eqref{thm2fast} in parts as two separate lemmas. The first concerns the linear part, and has been adapted from Escobedo-Zuazua~\cite{Escobedo-Zuazua} to apply to the fractional heat kernel.

\begin{lem} \label{EZ3} (\cite{Escobedo-Zuazua}) 
Let $p \in [1,\infty]$, and $\al \in [1,2]$. 
Let $\phi, |x|\phi \in L^1$, with $M := \intr \phi(x) \dd x$. 
Then there exists $C>0$, such that
\begin{align}
    \norm{\P(t) \ast \phi - M\P(t)}_p \leq C \norm{\phi}_{L^1(\R^2 ; |x|)} t^{\frac{2}{\al} (1 - \frac{1}{p}) - \frac{1}{\al} } . 
\end{align}
\end{lem}

\begin{lem} \label{nonlin est}
Let $ 1 \leq p \leq \infty$, and $\al\in[1,2]$. 
Suppose that $\theta$ is a mild solution satisfying the decay properties 
in Proposition~\ref{prop:decay}. Then 
there exists $C >0$ such that, for all $t\geq1$,
\begin{align*}
\bignorm{ \nonlin }_p  \leq Cb_{\alpha ,p} (t), 
\end{align*}
where $b_{\alpha ,p}$ is defined by \eqref{b}. 
\end{lem}

\begin{proof}
\begin{step}
($p \geq 2$ case)  
We will need to split the time interval into two halves, and handle the $\al = 1$ and $\al > 1$ cases separately. 
The proof below is only for the $p < \infty$ case, as the $p = \infty $ case is almost identical. 
\end{step}

We utilise the $L^p$ decay of the solution 
and its derivative. 
For values $t \in (0,1]$, we can increase the powers of $s$ in our estimates by taking 
\[
\norm{\t(t)}_{p} \leq \norm{\t_0}_{p} \leq C, \text{ for all } t>0,
\]
by which we ensure that the time-integral does not blow up locally.

We start on the second half of the time-interval, with $\al>1.$ 
By the boundedness of the Riesz transform, 
\begin{align*}
    \bignorm{\int^t_{t/2} \P(t-s) \ast \N \cdot (u(s)\t(s))  \dd s  }_p & \leq  \int^{t}_{t/2} \norm{\N \P(t-s)}_{1} \norm{u(s)\t(s)}_{p}  \dd s \\
    & \leq \int^{t}_{t/2} C (t-s)^{-1/\al} \norm{\t(s)}_{2p}^2  \dd s \\
    & \leq Ct^{-1/\al +1 } t^{-\frac{4}{\al} (1-\frac{1}{2p})  } \\
    & = C t^{-\frac{2}{\al} (1- \frac{1}{p}) -\frac{3}{\al} +1 } .
\end{align*}
In the $\al = 1$ case,
\begin{align*}
    \bignorm{\int^t_{t/2} P_{1/2}(t-s) \ast \N \cdot (u(s)\t(s))  \dd s}_p & \leq \int^t_{t/2} \norm{P_{1/2}(t-s)}_1 \norm{ (u(s)\cdot\N)\t(s) }_p   \dd s \\
    & \leq \int^{t}_{t/2} C \norm{\t(s)}_{2p} \norm{\N \t(s)}_{2p}  \dd s \\
    & \leq Ct^{-2(1 - \frac{1}{p}) -2 },
\end{align*}
where we have used \eqref{Iwab Lp deriv decay} to handle the derivative of $\t$.

For the first half of the time-interval, we will distinguish between the $\al>1$ and $\al=1$ cases when it becomes necessary. We will take the Fourier transform inside the norm, and manipulate the resulting multipliers from the derivative and Riesz transform. 
\begin{align*}
    & \bignorm{\int^{t/2}_0 \P(t-s) \ast \N \cdot(u(s)\t(s))   \dd s }_p \\
    & = \bignorm{\int^{t/2}_{0} \sum_{j =1}^{2} \F^{-1} \Big{[} e^{-(t-s) |\x|^\al } \intr \x_j \frac{(-1)^j}{2} \Big{(} \frac{\eta_{3-j}}{|\eta|} + \frac{\x_{3-j} - \eta_{3-j}}{|\x - \eta|} \Big{)}   \th(s, \x - \eta) \th(s,\eta) \dd \eta \Big{]} \dd s }_p \\
    & = \bignorm{ \int^{t/2}_{0} \sum_{j = 1}^2 \frac{1}{2} \F^{-1} \Bigf{ e^{-(t-s)|\x|^\al} \intr  \x_j \eta_{3-j}  \frac{|\x-\eta|^2 - |\eta|^2 }{|\eta| |\x-\eta| (|\x-\eta| + |\eta| )}   \th(s, \x - \eta) \th(s,\eta) \dd\eta  }   \dd s }_p. \\
\end{align*}

We start with the $p=2$ case. By the Plancherel theorem,
\begin{align*}
    & \bignorm{ \int^{t/2}_{0} \P(t-s) \ast \nabla \cdot \Bigc{ u(s)\t(s) } \dd s }_2 \\
    & \leq \bignorm{ \int^{t/2}_{0} |\x|^2 e^{-(t-s) |\x|^\al} \frac{1}{2} \intr 
    \frac{ |\eta| |\x-\eta| + |\eta|^2 }{|\eta| |\x-\eta| (|\x-\eta| + |\eta|)} 
    { \Big| \hat{\t}(s,\x-\eta) \hat{\t}(s,\eta) \Big| 
    }
     \dd\eta \dd s }_2.
\end{align*}
The large multiplier inside the $\eta$ integral is easily estimated from above by 
\[ \frac{2}{|\x - \eta|}. \]
Next, we take the Littlewood-Paley decomposition of both $\t$ functions: 
\begin{align*}
 & \hat{\t}(s,\x-\eta) = \sum_{k\in \Z} \hat{\phi}_k (\x-\eta) \hat{\t}(s,\x-\eta) =: \sum_{k \in \Z} \hat{\t}_k(s,\x-\eta) ,  \\
 & \hat{\t}(s,\eta) = \sum_{l\in \Z} \hat{\phi}_l (\eta) \hat{\t}(s,\eta) =: \sum_{l\in \Z} \hat{\t}_l (s,\eta).
\end{align*}
We also split the $L^2$ norm by H\"older's inequality.
\begin{align}
    & \bignorm{ \int^{t/2}_{0} |\x|^2 e^{-(t-s) |\x|^\al} \frac{1}{2} \intr 
    \frac{ |\eta| |\x-\eta| + |\eta|^2 }{|\eta| |\x-\eta| (|\x-\eta| + |\eta|)} \Big| \hat{\t}(s,\x-\eta) \hat{\t}(s,\eta)\Big|  \dd\eta \dd s }_2 \notag \\
    & \leq \sum_{k,l \in \Z} \int^{t/2}_{0} \norm{ |\xi|^2 e^{-(t-s)|\xi|^\alpha} }_2 \bignorm{  \intr  \frac{1}{|\x-\eta|}\Big|  \hat{\t}_k(s,\x-\eta) \hat{\t}_l (s,\eta) \Big|  \dd \eta  }_\infty  \dd s \notag \\
    & \leq Ct^{-\frac{2}{\al} (1-\frac{1}{2}) -\frac{2}{\al} } \int^{t/2}_{0} \sum_{k,l \in \Z} \frac{1}{2^k} \norm{\hat \t_k(s)}_{\frac{4}{3}} \norm{\hat \t_l (s)}_{4}  \dd s \notag \\
    & \leq Ct^{-\frac{2}{\al} (1-\frac{1}{2}) -\frac{2}{\al} } \int^{t/2}_{0} \sum_{k,l \in \Z} \norm{\t_k(s)}_{\frac{4}{3}} \norm{\t_l (s)}_{\frac{4}{3}}  \dd s. \label{prebesov}
\end{align}
The final step is obtained by applying the Hardy-Littlewood-Sobolev inequality.
The above sum can be written as the product of Besov norms of the solution, $\t$. We give estimates for the necessary Besov norms next.

Recall that $\t$ is made up of a linear and nonlinear term. We will take the norms of each separately, and will see that both hinge on the Besov norm of the fundamental solution, $\P$. We start with the linear term. Clearly
\begin{align*}
    \norm{\P(t)\ast\t_0}_{\dot{B}^{0}_{\frac{4}{3},1}} \leq  \norm{\P(t)}_{\dot{B}^{0}_{\frac{4}{3},1}}\norm{\t_0}_{1}
    \leq C t^{-\frac{1}{2\alpha}} \| \t_0 \|_1 . 
\end{align*}

Let us now consider the Besov norms of the whole solution, $\t$. In the following, we must distinguish between the noncritical and critical cases. We start with the noncritical $\al \in (1,2]$ case. 
\begin{align*}
    & \norm{\t(t)}_{\dot{B}^{0}_{\frac{4}{3},1}} 
    \leq C t^{-\frac{1}{2\alpha}} + \bignorm{\nonlin}_{\dot{B}^{0}_{\frac{4}{3},1}} \\
    & \leq C t^{-\frac{1}{2\alpha}} + \int^{t/2}_{0} \norm{\nabla\P(t-s)}_{\dot{B}^{0}_{\frac{4}{3},1}} \norm{u(s)\t(s)}_{1}  \dd s  \\
    & \ \ \ \ \ \ \ \ \ \ + \int^{t}_{t/2} \norm{\nabla\P(t-s)}_{\dot{B}^{0}_{1,1}} \norm{u(s)\t(s)}_{\frac{4}{3}}  \dd s \\
    & \leq Ct^{-\frac{1}{2\al} } + \int^{t/2}_{0} C(t-s)^{-\frac{1}{2\al} - \frac{1}{\al} } \norm{u(s)\t(s)}_{1}  \dd s + \int^{t}_{t/2} C(t-s)^{-1/\al} \norm{u(s)\t(s)}_{\frac{4}{3}}  \dd s \\
    & \leq Ct^{-\frac{1}{2\al} } + Ct^{-\frac{1}{2\al}  - \frac{1}{\al} } \int^{t/2}_{0} (s+1)^{-2/\al} \dd s + Ct^{-\frac{1}{2\al}  - \frac{2}{\al} } \int^{t}_{t/2} (t-s)^{-1/\al} \dd s \\
    & \leq Ct^{-\frac{1}{2\al}}, \text{ for all } t\geq 1, \text{ and all } \al \in (1,2].
\end{align*}
We now look at the critical case, $\al=1$. Here, the only difference is that 
we leave the derivative on the right hand side of the convolution 
for the second half of the time-integral, 
and apply the decay estimate of the derivative. 
\begin{align*}
    & \| \theta(t) \|_{\dot B^0_{\frac{4}{3},1} } 
     \leq Ct^{-\frac{1}{2}}  + \int^{t/2}_{0} \norm{\nabla P_{1/2}(t-s)}_{\dot{B}^{0}_{\frac{4}{3},1}} \norm{u(s)\t(s)}_{1}  \dd s  \\
    & \ \ \ \ \ \ \ \ \ \ \ \ \ \ \ \ + \int^{t}_{t/2} \norm{P_{1/2}(t-s)}_{1} \norm{u(s)\cdot\nabla\t(s)}_{\frac{4}{3}}  \dd s \\
    & \leq Ct^{-\frac{1}{2}} + \int^{t/2}_{0} C(t-s)^{-\frac{3}{2} } \norm{u(s)\t(s)}_{1}  \dd s  + \int^{t}_{t/2} C \norm{\t(s)}_{\frac{8}{3}}\norm{\nabla\t(s)}_{\frac{8}{3}}  \dd s \\
    & \leq Ct^{-\frac{1}{2}} + Ct^{-\frac{3}{2} } \int^{t/2}_{0} (s+1)^{-2} \dd s + C\int_{t/2}^t s^{-\frac{7}{2}} \dd s \\
    & \leq Ct^{-\frac{1}{2} }, \text{ for all } t\geq 1.
\end{align*}

Returning to \eqref{prebesov}, we get
\begin{align*}
 t^{-\frac{2}{\al} (1-\frac{1}{2}) -\frac{2}{\al} } \int^{t/2}_{0} \| \theta(s)\|_{\dot B^0_{\frac{4}{3},1}} ^2 
    \dd s 
    & \leq Ct^{-\frac{2}{\al} (1-\frac{1}{2}) -\frac{2}{\al} }
    \Big\{ \int^1_0 \dd s + \int^{t/2}_{1} s^{-\frac{1}{\alpha}} \dd s 
    \Big\} \\
    & \leq 
    %\begin{cases}
    %Ct^{1-\frac{4}{\al}} , \text{ for } \al \in (1,2], \\
    %Ct^{-3} \ln(t), \text{ for } \al = 1.
    %\end{cases}
    1 / b_{\al,2}(t).
\end{align*}

The above result is easily extended to all $p>2$ by Young's convolution inequality. We simply split the fundamental solution into two parts as follows, and then proceed through the exact same steps as above.
\begin{align*}
    & \bignorm{ \int^{t/2}_{0} \P(t-s) \ast \nabla \cdot \Bigc{ u(s)\t(s) } \dd s }_p \\
    & \leq \int^{t/2}_{0} \norm{\P((t-s)/2)}_{\frac{2p}{2+p}} \\
    & \ \ \ \ \bignorm{e^{-((t-s)/2)|\xi|^\al} \frac{1}{2} \intr (\x_2\eta_1 - \x_1\eta_2) \frac{ \x \cdot (\x-2\eta ) }{|\eta| |\x-\eta| (|\x-\eta| + |\eta|)} \hat{\t}(s,\x-\eta) \hat{\t}(s,\eta) \dd \eta}_2  \dd s.
\end{align*}

\begin{step}
($1 \leq p < 2$ case) 
We split the Fourier multiplier into two parts:
\[ \frac{ \eta_2 (\x_1 - 2\eta_1) }{|\eta| |\x - \eta| (|\x - \eta| + |\eta|)} = 
\frac{ \eta_2 (\x_1 - \eta_1) + \eta_1 \eta_2 }{|\eta| |\x - \eta| (|\x - \eta| + |\eta|)},\]
and we focus on the first term above, as both cases have almost identical proofs. We write 
\[
m(\xi-\eta,\eta) := 
\frac{ \eta_2 (\x_1 - \eta_1)  }{|\eta| |\x - \eta| (|\x - \eta| + |\eta|)} .
\]
Next, the most crucial step to this method is to split the Euclidean space into squares, whose size depends on time. For $k = (k_1,k_2) \in \mathbb Z^2$, define the set 
\begin{align} \notag 
Q_{t,k} \coloneqq \big{\{} x=(x_1,x_2) \in \R^2 \ | \ &   \ x_j \in [t^{1/\al}k_j, t^{1/\al}(k_j + 1)), \ j = 1,2 
 \big{\}}.
\end{align}
Then we can once again fit our inverse Fourier transform into an $L^2$ norm, and thus estimate away the Riesz transforms. That is,
\begin{align}
    & \bignorm{  \F^{-1} \Big{[} \int^{t/2}_{0} e^{-(t-s)|\x|^\al}  \intr \x^2_1 \, m(\xi-\eta,\eta ) \hat{\t}(s,\x-\eta) \hat{\t}(s,\eta)  \dd\eta \dd s \Big{]}  }_p \label{lessthan2} \\
    & \leq \sum_{k\in\Z^2} \bignorm{  \F^{-1} \Big{[} \int^{t/2}_{0} e^{-(t-s)|\x|^\al}  \intr \x^2_1 \, m(\xi-\eta,\eta ) \hat{\t}(s,\x-\eta) \hat{\t}(s,\eta)   \dd\eta \dd s \Big{]} }_{L^p(Q_{t,k})} \notag \\
    & \leq \sum_{k,l\in\Z^2} t^{\frac{2}{\al}(\frac{1}{p}-\frac{1}{2})} \bignorm{  \F^{-1} \Big{[} \int^{t/2}_{0} e^{-(t-s)|\x|^\al}  \intr  \xi_1 ^2 \, m(\xi-\eta,\eta ) \widehat{\t_{Q_{t,l}}}(s,\x-\eta) \hat{\t}(s,\eta)  
      \dd\eta \dd s \Big{]} }_{L^2(Q_{t,k})} . \notag 
\end{align}
where we have used H\"older's inequality. 
We are able to eliminate the Riesz transform by taking advantage of the Plancherel theorem, and split one of the solution functions, $\t$, into parts defined on squares, $Q_{t,l}$, as follows:
\[ \t = \sum_{l \in \Z^2} 1_{Q_{t,l}} \t =: \sum_{l\in\Z^2} \t_{Q_{t,l}}.
\] 
Next, we split up the above double sum into two cases: $k=l,$ and $k\neq l.$ 

In the case when $k=l$, 
\begin{align}
    & \sum_{ k \in \Z^2} t^{\frac{2}{\al} (\frac{1}{p}-\frac{1}{2})}  \bignorm{  \F^{-1} \Big{[} \int^{t/2}_{0} e^{-(t-s)|\x|^\al}  \intr 
\xi_1 ^2 \, m(\xi-\eta,\eta ) \widehat{\t_{Q_{t,k}}}(s,\x-\eta) \hat{\t}(s,\eta)     \dd\eta \dd s \Big{]} }_{L^2(Q_{t,k})} \notag \\
    & \leq \sum_{ k \in \Z^2}  t^{\frac{2}{\al}(\frac{1}{p}-\frac{1}{2})}  \int^{t/2}_{0} \norm{e^{-(t-s)|\x|^\al} \xi_1 ^2}_2 \bignorm{ \intr 
    \frac{ \Big| \widehat{\t_{Q_{t,k}}}(s,\x-\eta) \hat{\t}(s,\eta)
           \Big| }{ |\x - \eta|^{1/2}  |\eta|^{1/2}}   \dd\eta}_\infty \dd s \notag \\
    & \leq C t^{\frac{2}{\al p} - \frac{4}{\al} } \int^{t/2}_{0} \sum_{k\in\Z^2} \norm{ |\nabla|^{-1/2} \t_{Q_{t,k}}(s) }_{2} \norm{ |\nabla|^{-1/2} \t(s) }_2    \dd s \notag \\
    & \leq C t^{\frac{2}{\al p} - \frac{4}{\al} } \int^{t/2}_{0} \sum_{k\in\Z^2} \norm{ \t_{Q_{t,k}}(s) }_{4/3} \norm{ \t(s) }_{4/3}  \dd s. \label{how now brown cow} 
\end{align}
The steps above have used Young's convolution inequality and the Hardy-Littlewood-Sobolev inequality, as we have seen before in our original estimations of the nonlinear term of $\t$. 

In the $k \not = l$ case, we multiply by $1$ by inserting 
$|t^{1/\al}k - t^{1/\al}l|^{2} \text{/} |t^{1/\al}k - t^{1/\al}l|^{2}$, and thus write 
\begin{align}
    &  \sum_{ l \in \Z^2} \sum_{ k\neq l} t^{\frac{2}{\al} (\frac{1}{p}-\frac{1}{2})}   \frac{ |t^{1/\al}k - t^{1/\al}l|^{2} }{  |t^{1/\al}k-t^{1/\al}l|^{2}} \notag \\
    & \ \ \ \ \ \ \ \ \bignorm{  \F^{-1} \Big{[} \int^{t/2}_{0} e^{-(t-s)|\x|^\al}  \intr \xi_1 ^2 \, m(\xi-\eta,\eta ) \widehat{\t_{Q_{t,l}}}(s,\x-\eta) \hat{\t}(s,\eta) \dd\eta \dd s \Big{]} }_{L^2(Q_{t,k})} \notag \\
    & \leq \sum_{ l \in \Z^2} C t^{\frac{2}{\al} (\frac{1}{p}-\frac{1}{2}) - \frac{2}{\al} } \Biggc{ \sum_{ k\neq l} \Bigg{\{} 
    | t^{1/\al} k- t^{1/\al} l|^{2} \label{errorlol} \\
    & \ \ \ \ \ \ \ \ \bignorm{  \F^{-1} \Big{[} \int^{t/2}_{0} e^{-(t-s)|\x|^\al}  \intr \xi_1 ^2 \, m(\xi-\eta,\eta ) \widehat{\t_{Q_{t,l}}}(s,\x-\eta) \hat{\t}(s,\eta) \dd\eta \dd s \Big{]} }_{L^2(Q_{t,k})}
    \Bigg{\}}^2  }^{1/2} \notag
\\
    & \leq \sum_{ l \in \Z^2} C t^{\frac{2}{\al} (\frac{1}{p}-\frac{1}{2}) - \frac{2}{\al} } 
      \label{omg} \\
    & \ \ \ \ \ \ \ \ \bignorm{ | x- t^{1/\al} l|^{2} \F^{-1} \Big{[} \int^{t/2}_{0} e^{-(t-s)|\x|^\al}  \intr \xi_1 ^2 \, m(\xi-\eta,\eta ) \widehat{\t_{Q_{t,l}}}(s,\x-\eta) \hat{\t}(s,\eta) \dd\eta \dd s \Big{]} }_{2} . \notag 
\end{align}
The step \eqref{errorlol} was obtained by simply using H\"older's inequality for sequences, noting that
\[
 \sum_{k\neq l} |t^{1/\al}k - t^{1/\al}l|^{-4} \leq C t^{-4/\al}, \text{ for all } \al\in[1,2], \ t\geq1, \ k,l\in\Z^2.
\]
Our next concern is with the boundedness of the sum over $l.$ The key points to the following steps are that $t^{1/\al} k$ is close to $x$, where $x$ is the variable of our $L^2(Q_{t,k})$ norm; and that we thus treat $t^{1/\al} k$ as a derivative in $\x$ after moving it inside the inverse Fourier transform. For each $l$, 
\begin{align*}
    & 
\bignorm{ | x- t^{1/\al} l|^{2} \F^{-1} \Big{[} \int^{t/2}_{0} e^{-(t-s)|\x|^\al}  \intr \xi_1 ^2 \, m(\xi-\eta,\eta ) \widehat{\t_{Q_{t,l}}}(s,\x-\eta) \hat{\t}(s,\eta) \dd\eta \dd s \Big{]} }_{2}  \\
    & =\bignorm{  \F^{-1} \Big{[} ( i \nabla_\x - t^{1/\alpha}l )^{2} \int^{t/2}_{0} e^{-(t-s)|\x|^\al}  \intr \xi_1 ^2 \, m(\xi-\eta,\eta ) \widehat{\t_{Q_{t,l}}}(s,\x-\eta) \hat{\t}(s,\eta) \dd\eta \dd s \Big{]} }_{2}.
\end{align*}
Briefly, the $\x$ derivative results in a factor of $t^{1/\al},$ which simply cancels with the inverse factor emerging from the $\x_1^2$ already present inside the integral. The resulting estimate is thus essentially similar to
\begin{align*}
    & 
    \bignorm{  \F^{-1} \Big{[} \int^{t/2}_{0} e^{-(t-s)|\x|^\al}  \intr  \, m(\xi-\eta,\eta ) \widehat{\t_{Q_{t,l}}}(s,\x-\eta) \hat{\t}(s,\eta) \dd\eta \dd s \Big{]} }_{2}
    \\
    & \leq C \int^{t/2}_0 \norm{\P(t-s)}_2 \bignorm{ \intr  \frac{ \Big|\widehat{\t_{Q_{t,l}}} (s,\x-\eta) \hat{\t}(s,\eta)\Big| }{( |\x-\eta| + |\eta| )} \dd \eta }_\infty  \dd s \\
    & \leq C t^{-1/\al} \int^{t/2}_{0} \norm{ |\nabla|^{-1/2} \t_{Q_{t,l}} (s) }_2 \norm{ |\nabla|^{-1/2} \t(s) }_2   \dd s \\
    & \leq C t^{-1/\al} \int^{t/2}_{0}  \norm{\t_{Q_{t,l}}(s) }_{4/3} \norm{\t(s)}_{4/3}   \dd s. 
\end{align*}
Returning to \eqref{omg}, 
\begin{align}
    &  \sum_{ l \in \Z^2} \sum_{ k\neq l} t^{\frac{2}{\al} (\frac{1}{p}-\frac{1}{2})}   \frac{ |t^{1/\al}k - t^{1/\al}l|^{2} }{  |t^{1/\al}k-t^{1/\al}l|^{2}} \notag \\
    & \ \ \ \ \ \ \ \ \bignorm{  \F^{-1} \Big{[} \int^{t/2}_{0} e^{-(t-s)|\x|^\al}  \intr \xi_1 ^2 \, m(\xi-\eta,\eta ) \widehat{\t_{Q_{t,l}}}(s,\x-\eta) \hat{\t}(s,\eta) \dd\eta \dd s \Big{]} }_{L^2(Q_{t,k})} \notag \\
& 
\leq 
 C t^{\frac{2}{\al p} - \frac{4}{\al} } \int^{t/2}_{0} \sum_{l\in\Z^2} \norm{ \t_{Q_{t,l}}(s) }_{4/3} \norm{ \t(s) }_{4/3}  \dd s .
 \label{330-1}
\end{align}

It is a delicate process to prove that the above sums \eqref{how now brown cow} and \eqref{330-1} are bounded properly. As such, we provide an outline of the proof, beginning with the following estimates.
\[
\begin{split}
\sum_{k\in\Z^2} \norm{  \t_{Q_{t,k}}(t) }_{4/3} 
\leq 
& 
C t^{-\frac{1}{2\alpha}} 
+ C t^{\frac{2}{\alpha p}-\frac{4}{\alpha}} 
\int^{t/2}_{0} \sum_{k\in\Z^2} \norm{ \t_{Q_{t,k}}(s) }_{4/3} \norm{ \t(s) }_{4/3}  \dd s
\\
& 
+ C \int_{t/2}^t \sum_{k\in\mathbb Z^2} \| 1_{Q_{t,l}}P_{\alpha/2}(t- s ) \|_{\frac{4}{3}} \| (u(s) \cdot \nabla) \theta(s) \|_{1} \dd s 
\\
\leq 
& 
C t^{-\frac{1}{2\alpha}} 
+ C t^{-\frac{1}{2\alpha}} 
\int^{t/2}_{0} \sum_{k\in\Z^2} \norm{ \t_{Q_{t,k}}(s) }_{4/3} (1+s)^{-\frac{3}{2\alpha}}   \dd s, \quad t \geq 1,
\\
\leq 
& 
C t^{-\frac{1}{2\alpha}} 
+ C t^{-\frac{1}{2\alpha}} 
\int^{t/2}_{0} \sum_{k\in\Z^2} \norm{ \t_{Q_{s,k}}(s) }_{4/3} (1+s)^{-\frac{3}{2\alpha}}   \dd s, \quad t \geq 1,
\end{split}
\]
where here we have applied the decay estimate of the solution and its derivative,
and the elementary inequality 
$t^{-1} \leq C (1+s)^{-1}$, for $t \geq \max\{ 1, s \}$. 
By Gr\"onwall's inequality we have the same decay as in $L^{\frac{4}{3}}$, 
\[
\sum_{k\in\Z^2} \norm{  \t_{Q_{t,k}}(t) }_{4/3} 
\leq C t^{-\frac{1}{2\alpha}} , \quad t \geq 1. 
\]
We apply the decay above to the inequalities \eqref{how now brown cow} and \eqref{330-1}, and so \eqref{lessthan2} 
is bounded by 
\[
C t^{\frac{2}{\alpha p}-\frac{4}{\alpha}} 
\int_0^{t/2} (1+s)^{-\frac{1}{2\al}} (1+s)^{-\frac{1}{2\alpha}} \dd s 
= 
 1 / b_{\al,p}(t).
\]\end{step} 
This completes the proof in the case when $1 \leq p < 2$. 
\end{proof}

We thus have completed the proof of \eqref{thm2fast}, and also obtained a useful bound on the nonlinear term of our solution, which we will use later. For all $p \in [1,\infty]$, and with initial data satisfying \eqref{data}, 
\begin{align} \label{nonlin bound}
    \bignorm{\nonlin}_p \leq 
    1 / b_{\al,p}(t)
\end{align}
for all $\al \in [1,2], \ t \geq 1$.

\subsection{Improving the Decay Rate}

We now begin our proof of Theorem 1.1. 
We will show convergence for the linear and nonlinear approximations separately. The convergence of the linear part is expressed in the next lemma. 
\begin{lem}
Let $\al \in [1,2]$ and $p \in [1,\infty]$. Let $\phi, |x|^2\phi \in L^1$ and $M := \intr \phi(x) \dd x$. Then
\begin{align*}
    t^{\frac{2}{\al} (1-\frac{1}{p}) + \frac{2}{\al} } 
    & \norm{ \phi \ast \P(t) - M\P(t) + \nabla \P(t) \cdot \intr y \phi(y) \dd y \\
    & - \sum_{i,j = 1}^2 \partial_{x_i} \partial_{x_j} \P(t) \intr y_i y_j \phi(y) \dd y }_{p} \to 0, \text{ as } t \to \infty.
\end{align*}
\end{lem}

\begin{proof}
This lemma is proven analogously to Lemma~\ref{EZ3}. Indeed, 
by the Taylor expansion 
\[
P_{\frac{\alpha}{2}} (t, x-y) 
= P_{\alpha/2}(t,x) 
-\nabla P_{\alpha /2}(x) \cdot y 
+ \sum_{i,j = 1}^2 \partial_{x_i} \partial_{x_j} \P(t) y_i y_j ,
\quad y \to 0 ,
\]
the convergence result becomes clear; and so we omit the details. 
\end{proof}

Finally, we discuss the nonlinear approximation.
\begin{lem}
Let $\al \in [1,2]$ and $p \in [1,\infty]$. 
\begin{align*}
    &U(t) := \lin, \quad  
    I(t) := \nonlin.
\end{align*}
Then
\begin{align*}
    b_{\alpha ,p}(t) 
     \bignorm{ I(t) - \int^{t}_{0} \P(t-s) \ast \nabla \cdot \Bigc{(RU(s)) U(s)} \dd s  }_{p} \to 0, \text{ as } t \to \infty. \end{align*}
\end{lem}

\begin{proof}

The convergence can be shown by a similar method to the bounds from above that we have calculated up to this point. We split the time interval into two halves. The second half is estimated simply using H\"older's inequality and Young's convolution inequality. The first half is estimated using the same method as in Section 3.1.

The key point is that a faster decay is achieved when taking the difference between the two terms above than when they are estimated separately. 
This is accomplished by splitting the difference as follows:
\begin{align}
    & 
    (u(s)\cdot \nabla )\t(s) - \Bigc{RU(s)\cdot \nabla}U(s) 
= N\big(U(s) , -I(s)\big) - N\big(I(s) , \t (s) \big), 
\label{two parts}
\end{align}
where $R := (-R_2, R_1)$, and 
\[
N(f,g) := 
\sum_{j \neq k} \F^{-1} \Big{[} \intr \x_j \frac{(-1)^j}{2} \Big{(} \frac{\eta_k}{|\eta|} + \frac{\x_{k} - \eta_k}{|\x - \eta|} \Big{)}   \hat{f}( \x - \eta) \hat{g}(\eta) \dd \eta
\Big{]}.
\]
We know that there are faster decay estimates for the two terms in the right hand side of \eqref{two parts}
than for $(u \cdot \nabla) \theta$, since the nonlinear part 
$I$ decays faster than the linear part.

We split the time interval into two halves again. The second half is simply calculated using estimates of $\t,U,I$ that we have seen above.
For the first half of the time interval, we again split up the proof into the $p\geq2$ and $1\leq p <2$ cases.

\setcounter{step}{0}

\begin{step}
(First half of time interval, $p\geq2$ case) We will show the proof for $p=2,$ and again the result can be easily extended to all greater values of $p$ afterwards. We handle the two terms in \eqref{two parts} separately.
\begin{align*}
    & \bignorm{ \int^{t/2}_{0} |\xi|^2 e^{-(t-s) |\x|^\al} \intr \frac{ |\eta| |\x-\eta| + |\eta|^2 }{|\eta| |\x-\eta| (|\x-\eta| + |\eta|)} \hat{U}(s,\x-\eta) \hat{I}(s,\eta)   \dd \eta   \dd s }_2  \\
    & \leq C t^{-\frac{2}{\al} (1-\frac{1}{2}) - \frac{2}{\al} } \int^{t/2}_{0} \norm{U (s)}_{\dot B^0_{\frac{4}{3},1}} \norm{I (s)}_{\dot B^0_{\frac{4}{3},1}} \dd s \\
    & \leq Ct^{-3/\al} 
    \Big(\int_0^1 \dd s +  \int^{t/2}_{1} (1+s)^{-\frac{1}{2\alpha}} \cdot b_{\alpha,\frac{4}{3}}(s)\dd s\Big),
\end{align*} \end{step}
\noindent which proves that 
\[
b_{\alpha ,2} (t)\bignorm{ \int^{t/2}_{0} P_{\alpha /2}(t) *  N\big(U(s) , -I(s)\big)
 \dd s }_2 
 \to 0 , \quad t\to\infty .
\]
The convergence of the second term 
$N\big( I(s), \theta(s) \big)$ follows from almost the same argument as above, 
by applying the decay of $\theta$ instead of $U$.

\begin{step}
(First half of time interval, $1 \leq p < 2$ case) 
We again take the terms from \eqref{two parts}, 
and use the same technique as was used for the bound from above in section 3.1. 
That is, we split the Fourier multiplier 
and divide the time-integral into two halves, 
and split the space into squares $Q_{t,k}$. 
We omit the details. 
\end{step}\end{proof}

\subsection{Optimal Decay of NonLinear Term}

We lastly discuss the optimality (in the $p=2$ case) of our estimate in Lemma \ref{nonlin est}. The decay rate for our estimate from above is optimal if we can bound the nonlinear estimate from below by the same power of $t.$ That is, we need
\begin{align} 
    & \bignorm{\intt \P(t-s) \ast \N \cdot \Big{(} (RU(s))U(s) \Big{)} \dd s}_2 \geq 
    b_{\al,2}(t),
\end{align}
for all $\al \in [1,2]$, and $t>0$ sufficiently large.

Since we are taking the $L^2$-norm, taking the Fourier Transform inside the norm does not change its value.
\begin{align*}
& \bignorm{ \int^{t}_{0} P_{\alpha/2}(t-s) \ast \nabla \cdot ((RU(s))U(s)) \dd s }_2 \\
&=
\bignorm{ \int^{t}_{0} \sum_{j = 1}^2 \xi_j e^{-(t-s)|\xi|^\alpha} \frac{(-1)^j}{2} \int_{\mathbb{R}^2} \Bigg{(} \frac{\eta_{3-j}}{|\eta|} + \frac{\xi_{3-j} - \eta_{3-j}}{|\xi - \eta|} \Bigg{)} e^{-s|\xi - \eta|^\alpha} e^{-s|\eta|^\alpha} \\
& \quad \hat{\theta}_0(\xi - \eta) \hat{\theta}_0(\eta) \dd \eta \dd s  }_2 .
\end{align*}

We rewrite the divergence operator and Riesz transform as two separate Fourier multipliers. 
\begin{align*}
& 
 \sum_{j = 1}^2 \xi_j  \frac{(-1)^j}{2} \Bigg{(} \frac{\eta_{3-j}}{|\eta|} + \frac{\xi_{3-j} - \eta_{3-j}}{|\xi - \eta|} \Bigg{)}
 \\
= &
 \frac{2\xi_1\xi_2}{|\xi - \eta|}
\bigg{(}\frac{\eta_1^2 - \eta_2^2}{|\eta|(|\eta|+|\xi - \eta|)}\bigg{)}
 + \frac{2\eta_1\eta_2}{|\xi - \eta|}
\bigg{(}\frac{\xi_2^2 - \xi_1^2}{|\eta|(|\eta|+|\xi - \eta|)}\bigg{)}   \\
& + \sum_{j =1}^2 \frac{|\xi|^2 \xi_j \eta_{3-j}}{|\eta||\xi-\eta|(|\eta|+|\xi-\eta|)}
\\
=: 
& m_1 (\xi-\eta ,\eta) + m_2 (\xi-\eta , \eta)
\end{align*}

The key difference between these two multipliers is that the numerator of $m_1$ features a second-order derivative, whereas that of $m_2$ has a third-order derivative.
We show that, for some initial data $\theta _0$, the first part with $m_1$ has the optimal decay and the remainder with $m_2$ is smaller. 

\begin{lem}
Let $\delta, \epsilon > 0$. Let $\t_0 \in W^{1,1} \cap W^{1,\infty}$ as before, but with the following additional conditions: 
\begin{itemize}
    \item $\th_0 \geq 0, \text{ on } \R^2$, 
    \item $\supp \th_0 \subseteq \{ \x \in \R^2 \ | \ |\x_2| < \delta |\x_1| \}$, 
    \item $\th_0(\x) \geq C, \text{ for some } C>0, \text{ for all } \xi \in \supp \th \cap \{ \x \in \R^2 \ | \ |\xi| \leq 1 \}$.
\end{itemize}
Then, for sufficiently small $\delta$ and $\epsilon$, we have 
\begin{align} 
\label{m_1 larger}
\bignorm{ \int^{t}_{0} e^{-(t-s)|\xi|^\alpha}  
\int_{\mathbb{R}^2} m_1 (\xi-\eta,\eta) \, e^{-s|\xi - \eta|^\alpha} e^{-s|\eta|^\alpha} \dd \eta \dd s  }_{L^2(|\xi| \leq \epsilon t^{-1/\alpha})}
\notag \\
 \geq \begin{cases}
C t^{1-\frac{4}{\al}} \epsilon^{3}, \text{ for } \al \in (1,2], \\
C t^{-3} \epsilon^{3} \ln(t) , \text{ for } \al =1, 
\end{cases}
\\
\label{m_2 smaller}
\bignorm{ \int^{t}_{0} e^{-(t-s)|\xi|^\alpha}  
\int_{\mathbb{R}^2} m_2 (\xi-\eta,\eta) \, e^{-s|\xi - \eta|^\alpha} e^{-s|\eta|^\alpha} \dd \eta \dd s  }_{L^2(|\xi| \leq \epsilon t^{-1/\alpha})} 
\notag \\
\leq C \epsilon^{7/2} t^{1-\frac{4}{\alpha}}, \text{ for all } \al \in [1,2].
\end{align}

\end{lem}

\begin{proof}
We consider \eqref{m_2 smaller} and \eqref{m_1 larger} separately. Beginning with \eqref{m_2 smaller}, we consider just the $j=1$ part, as the estimates of both terms are identical.
\begin{align*}
& \bignorm{ \int^{t}_{0} e^{-(t-s)|\xi|^\alpha}  
\int_{\mathbb{R}^2} \frac{|\xi|^2 \xi_1 \eta_{2}}{|\eta||\xi-\eta|(|\eta|+|\xi-\eta|)} \hat{U}(s,\xi-\eta)\hat{U}(s,\eta) \dd \eta \dd s  }_{L^2(|\xi| \leq \epsilon t^{-1/\alpha})} \\
& \leq
\bignorm{ \int^{t}_{0} |\xi|^3 e^{-(t-s)|\xi|^\alpha}  
\int_{\mathbb{R}^2} \frac{ \eta_{2}}{|\eta||\xi-\eta|(|\eta|+|\xi-\eta|)} e^{-s|\xi - \eta|^\alpha} e^{-s|\eta|^\alpha} \hat{\theta}_0(\xi - \eta) \hat{\theta}_0(\eta) \dd \eta \dd s  }_{L^2(|\xi| \leq \epsilon t^{-1/\alpha})} \\
& \leq C 
\bignorm{ \int^{t}_{0} |\xi|^3 e^{-(t-s)|\xi|^\alpha}  
\int_{\mathbb{R}^2} \frac{ \eta_{2}}{|\eta||\xi-\eta|(|\eta|+|\xi-\eta|)} e^{-s|\xi - \eta|^\alpha} e^{-s|\eta|^\alpha} \dd \eta \dd s  }_{L^2(|\xi| \leq \epsilon t^{-1/\alpha})},
\end{align*}
as $\hat{\theta}_0$ is bounded. We bound the above norm by considering the integral
\[
\int_{\mathbb{R}^2} \frac{ \eta_{2}}{|\eta||\xi-\eta|(|\eta|+|\xi-\eta|)} e^{-s|\xi - \eta|^\alpha} e^{-s|\eta|^\alpha} \dd \eta.
\] 
Note that, by making the substitution $\eta \rightarrow \xi - \eta$, we can rewrite the integral as 
\begin{align*}
& \int_{\mathbb{R}^2} \frac{ \eta_{2}}{|\eta||\xi-\eta|(|\eta|+|\xi-\eta|)} e^{-s|\xi - \eta|^\alpha} e^{-s|\eta|^\alpha} \dd \eta \\
&= 
\frac{1}{2} \int_{\mathbb{R}^2} \frac{ (\eta_2 + (\xi_2 - \eta_{2}))}{|\eta||\xi-\eta|(|\eta|+|\xi-\eta|)} e^{-s|\xi - \eta|^\alpha} e^{-s|\eta|^\alpha} \dd \eta \\
&=
\frac{1}{2} \xi_2 
\Big( \int_{|\eta| \leq \frac{1}{2}|\xi|} +
\int_{\frac{1}{2}|\xi| \leq |\eta| \leq 2|\xi|} +
\int_{2|\xi| \leq |\eta|}
\Big) 
\frac{1}{|\eta||\xi-\eta|(|\eta|+|\xi-\eta|)} e^{-s|\xi - \eta|^\alpha} e^{-s|\eta|^\alpha} \dd \eta \\
&=: A_1 + A_2 + A_3,
\end{align*}
where we have split the integral into three parts with $|\eta|$ small, $|\eta|$ close to $|\x|$, and $|\eta|$ large. We start with the small part.
\begin{align*}
|A_1|
& \leq
C |\xi| \int_{|\eta| \leq \frac{1}{2}|\xi|} \frac{1}{|\xi|^2 |\eta|} \dd \eta \leq C . 
\end{align*}
Next we take $|\eta|$ large.
\begin{align*}
A_3
& \leq
C |\xi| \int_{|\eta| \geq 2|\xi|}\frac{1}{|\eta|^{3}} \dd \eta 
\leq C,
\end{align*}
and finally we take $|\eta|$ close to $|\x|$. 
\begin{align*}
A_2
& \leq 
C|\xi| \int_{\frac{1}{2}|\xi| \leq |\eta| \leq 2|\xi|}
\frac{e^{-s|\xi - \eta|^\alpha}}{|\xi|^2|\xi - \eta|} \dd \eta 
\leq 
C|\xi|^{-1} \int_{|\tilde \eta| \leq |\xi|} \frac{e^{-s|\tilde \eta|^\alpha}}{|\tilde \eta|}  \dd \tilde \eta \\
&
= C |\xi|^{-1} s^{-\frac{1}{\alpha}} \int_{|\tilde \eta| \leq s^{\frac{1}{\alpha}}|\xi|} \frac{e^{-s|\tilde \eta|^\alpha}}{|\tilde \eta|}  \dd \tilde \eta
 \leq 
C |\xi|^{-1} s^{-\frac{1}{\alpha}} \min \{ s^{\frac{1}{\alpha}}|\xi| ,1 \}
\\
& \leq C |\xi|^{-1/2} s^{-1/2\alpha} .
\end{align*}
Therefore, we obtain
\begin{align*}
&\bignorm{ \int^{t}_{0} |\xi|^3 e^{-(t-s)|\xi|^\alpha}  \int_{\mathbb{R}^2} \frac{ \eta_{2}}{|\eta||\xi-\eta|(|\eta|+|\xi-\eta|)} e^{-s|\xi - \eta|^\alpha} e^{-s|\eta|^\alpha} \dd \eta \dd s  }_{L^2(|\xi| \leq \epsilon t^{-1/\alpha})} \\
& \leq 
\bignorm{ \int^{t}_{0} |\xi|^3 e^{-(t-s)|\xi|^\alpha}  (C + C|\xi|^{-1/2}s^{-1/2\alpha})
\dd s  }_{L^2(|\xi| \leq \epsilon t^{-1/\alpha})} \\
& \leq
\bignorm{ |\xi|^3 \int^{t}_{0}  C \dd s  
}_{L^2(|\xi| \leq \epsilon t^{-1/\alpha})} 
+ 
\bignorm{ |\xi|^{5/2} \int^{t}_{0}  Cs^{-1/2\alpha}
\dd s  }_{L^2(|\xi| \leq \epsilon t^{-1/\alpha})} \\
& \leq
C(\epsilon^{4} + \epsilon^{7/2} )t^{1-\frac{4}{\alpha}} .
\end{align*}
We note that all of these terms are smaller than
\[
\epsilon^{7/2} t^{1-\frac{4}{\alpha}}
\]
for $\epsilon < 1$ and $t > 1$, and thus we have \eqref{m_2 smaller}. 

We next will show \eqref{m_1 larger}, for sufficiently small $\epsilon$, and for sufficiently large $t$.
We begin by labelling the two terms in our integral.

\newpage

\begin{align}& 
\bignorm{ \int^{t}_{0} e^{-(t-s)|\xi|^\alpha}  \nonumber
\int_{\mathbb{R}^2} 
\Bigg{(} \frac{2\xi_1\xi_2}{|\xi - \eta|}
\bigg{(}\frac{\eta_1^2 - \eta_2^2}{|\eta|(|\eta|+|\xi - \eta|)}\bigg{)}
 + \frac{2\eta_1\eta_2}{|\xi - \eta|}
\bigg{(}\frac{\xi_2^2 - \xi_1^2}{|\eta|(|\eta|+|\xi - \eta|)}\bigg{)} \Bigg{)} \\
&  \hat{U}(s,\xi-\eta)\hat{U}(s,\eta) \dd \eta \dd s  }_{L^2(|\xi| \leq \epsilon t^{-1/\alpha})} 
\label{wowee first term} \\
& =: \norm{(L) + (R)}_{L^2(|\xi| \leq \epsilon t^{-1/\al})}.
\end{align}
Our plan is to make $(L)$ the larger term. This is accomplished by the conditions on $\t_0$ that we have imposed. Taking the norm of $(L)$ on its own,
\begin{align*}
& \bignorm{ \int^{t}_{0} e^{-(t-s)|\xi|^\al}
\int_{\mathbb{R}^2} 
 \frac{2\xi_1\xi_2}{|\xi - \eta|}
\bigg{(}\frac{\eta_1^2 - \eta_2^2}{|\eta|(|\eta|+|\xi - \eta|)}\bigg{)} 
e^{-s|\x-\eta|^\al} \th_0(\x-\eta) e^{-s|\eta|^\al}\th_0(\eta) \dd \eta \dd s  }_{L^2(|\xi| \leq \epsilon t^{-1/\al})} \\
& \geq C \bignorm{ \xi_1\xi_2 \int^{t}_{0} e^{-(t-s)|\xi|^\al}
\int_{\mathbb{R}^2} 
\frac{\eta_1^2 e^{-s|\x-\eta|^\al} \th_0(\x-\eta) e^{-s|\eta|^\al}\th_0(\eta) }{|\eta| |\xi - \eta| (|\eta|+|\xi - \eta|)}
 \dd \eta \dd s  }_{L^2(|\xi| \leq \epsilon t^{-1/\al})} \\
 & \geq C \bignorm{ \xi_1\xi_2 \int^{t}_{1} e^{-\epsilon}
\int_{2 |\x| < |\eta| < 1} 
\frac{\eta_1^2 e^{-cs|\eta|^\al}  }{|\eta|^3}
 \dd \eta \dd s  }_{L^2(|\xi| \leq \epsilon t^{-1/\al})},
\end{align*}
where in the first step we made use of the shape of$\supp \th_0$, and in the second we used the fact that the integrands are positive and the bound from below for $\th_0$ on$\supp \th_0$ close to $0$. We next convert the integral over $\eta$ to polar coordinates and use substitution of variables to produce the final powers of $t$ in the subcritical case, and the $\ln$ function in the critical case.
\begin{align*}
& \bignorm{ \xi_1\xi_2 \int^{t}_{1} e^{-\epsilon}
\int_{2 |\x| < |\eta| < 1} 
\frac{\eta_1^2 e^{-cs|\eta|^\al}  }{|\eta|^3}
\dd \eta \dd s  }_{L^2(|\xi| \leq \epsilon t^{-1/\al})} \\
& \geq \bignorm{ \xi_1\xi_2 \int^{t}_{1} e^{-\epsilon}
\int^{1}_{2 |\x|} 
e^{-cs\rho^\al}
\dd \rho \dd s  }_{L^2(|\xi| \leq \epsilon t^{-1/\al})} \\
& = \bignorm{ \xi_1\xi_2 \int^{t}_{1} e^{-\epsilon}
\int^{s^{1/\al}}_{2 |\x| s^{1/\al}} 
e^{-c\rho}
\dd \rho s^{-1/\al} \dd s  }_{L^2(|\xi| \leq \epsilon t^{-1/\al})} \\
& \geq \bignorm{ \xi_1\xi_2 \int^{t}_{1} e^{-\epsilon}
C s^{-1/\al} \dd s  }_{L^2(|\xi| \leq \epsilon t^{-1/\al})} \\
& \geq \begin{cases}
C t^{1-\frac{4}{\al}} \epsilon^{3}, \text{ for } \al \in (1,2], \\
C  t^{-3} \epsilon^{3} \ln(t), \text{ for } \al =1.
\end{cases}
\end{align*}

Finally, by our setting of$\supp \th_0$, we obtain 
\begin{align*}
    \norm{(R)}_{L^2(|\xi| \leq \epsilon t^{-1})} \leq C  t^{-3} \epsilon^4 \ln{(t)},
\end{align*}
by estimations of integrals similar to before. Thus we obtain \eqref{m_1 larger}.
\end{proof}

\phantom{}

\noindent {\bf Conflict of interest statement.} On behalf of all authors, the corresponding author states that there is no conflict of interest. 

\phantom{}

\noindent {\bf Acknowledgements.} T. Iwabuchi was supported by JSPS Grant-in-Aid for Young Scientists
(A) (No.~17H04824).

\end{document}